\documentclass{article}
\usepackage[T1]{fontenc}    % Il faut utiliser le codage de fonte (font encoding)
\usepackage[english]{babel}
\usepackage{latexsym}
\usepackage{amssymb}
\usepackage{amsmath}
\usepackage{amsthm}
\usepackage{graphicx}
\usepackage{hyperref}
\usepackage{color}
\usepackage{colortbl}
\usepackage{epsfig}

\newtheorem{prop}{Proposition}[section] % Les propositions sont numérotées comme

\begin{document}

%============================== Title ================================

\title{SLOW INVARIANT MANIFOLD OF HEARTBEAT MODEL}
\author{{\it Jean-Marc GINOUX \& Bruno ROSSETTO}
\\[2ex] Universit\'{e} de Toulon,
\\ CNRS, LSIS, UMR 7296,
\\ B.P. 20132, 83957, La Garde Cedex France
\\ e-mail: {\tt ginoux@univ-tln.fr}
\\ {\tt http://ginoux.univ-tln.fr/}}

\maketitle

%============================ No page numbers =======================

\thispagestyle{empty}
\pagestyle{empty}

%============================ The Paper ==============================

%---------------------------- Abstract -------------------------------

\begin{abstract}
A new approach called \textit{Flow Curvature Method} has been recently developed in a book entitled \textit{Differential Geometry Applied to Dynamical Systems}. It consists in considering the trajectory curve, integral of any $n$-dimensional dynamical system as a curve in Euclidean $n$-space that enables to analytically compute the curvature of the trajectory - or the flow. Hence, it has been stated on the one hand that the location of the points where the curvature of the flow vanishes defines a manifold called \textit{flow curvature manifold} and on the other hand that such a manifold associated with any $n$-dimensional dynamical system directly provides its \textit{slow manifold analytical equation} the invariance of which has been proved according to Darboux theory. The \textit{Flow Curvature Method} has been already applied to many types of autonomous dynamical systems either singularly perturbed such as Van der Pol Model, FitzHugh-Nagumo Model, Chua's Model, ...) or non-singularly perturbed such as Pikovskii-Rabinovich-Trakhtengerts Model, Rikitake Model, Lorenz Model,... Moreover, it has been also applied to non-autonomous dynamical systems such as the Forced Van der Pol Model. In this article it will be used for the first time to analytically compute the \textit{slow invariant manifold analytical equation} of the four-dimensional Unforced and Forced Heartbeat Model. Its \textit{slow invariant manifold equation} which can be considered as a "state equation" linking all variables could then be used in heart prediction and control according to the strong correspondence between the model and the physiological cardiovascular system behavior.

\end{abstract}

%---------------------- General Instructions -------------------------

\section{Introduction}
\label{sec1} Dynamical systems consisting of \textit{nonlinear} differential
equations are generally not integrable. In his famous memoirs:
\textit{Sur les courbes d\'{e}finies par une \'{e}quation
diff\'{e}rentielle}, Poincar\'{e} \cite{poin1,poin4} faced to this problem
proposed to study \textit{trajectory curves} properties in the
\textit{phase space}.\\

``\ldots any differential equation can be written as:

\[
\frac{dx_1 }{dt} = X_1 ,
\quad
\frac{dx_2 }{dt} = X_2 , \quad \ldots ,
\quad
\frac{dx_n }{dt} = X_n
\]

where $X$ are integer polynomials.

If $t$ is considered as the time, these equations will define the motion of
a variable point in a space of dimension $n$.''

\begin{flushright}
-- Poincar\'{e} (1885, p. 168) --
\end{flushright}

Let's consider the following system of differential equations
defined in a compact E included in $\mathbb{R}$ as:

\begin{equation}
\label{eq1}
\frac{d\vec {X}}{dt} = \overrightarrow \Im ( \vec {X} )
\end{equation}

with

\[
\vec {X} = \left[ {x_1 ,x_2 ,...,x_n } \right]^t \in E \subset {\mathbb{R}}^n
\]

and

\[
\overrightarrow \Im ( \vec {X} ) = \left[ {f_1 ( \vec {X}),f_2 ( \vec {X} ),...,f_n ( \vec {X} )}
\right]^t \in E \subset {\mathbb{R}}^n
\]

The vector $\overrightarrow \Im ( \vec {X} )$ defines a velocity
vector field in E whose components $f_i $ which are supposed to be
continuous and infinitely differentiable with respect to all $x_i $ and $t$,
i.e. are $C^\infty $ functions in E and with values included in
$\mathbb{R}$, satisfy the assumptions of the Cauchy-Lipschitz theorem. For
more details, see for example \cite{cod}. A solution of this system
is a \textit{trajectory curve} $\vec {X}\left( t \right)$ tangent\footnote{ Except at the \textit{fixed points}.} to
$\overrightarrow \Im $ whose values define the \textit{states} of the \textit{dynamical system} described by the Eq.
(\ref{eq1}).

\newpage

Thus, \textit{trajectory curves} integral of dynamical systems
(\ref{eq1}) regarded as $n$-dimensional \textit{curves}, possess
local metrics properties, namely \textit{curvatures }which can be
analytically\footnote{ Since only time derivatives of the
\textit{trajectory curves} are involved in the \textit{curvature}
formulas. } deduced from the so-called Fr\'{e}net formulas \cite{fren}. For low dimensions two and three the concept of
\textit{curvatures} may be simply exemplified. A
three-dimensional\footnote{ A two-dimensional curve, i.e. a plane
\textit{curve} has a \textit{torsion} vanishing identically.} curve
for example has two \textit{curvatures: curvature} and
\textit{torsion} which are also known as \textit{first} and
\textit{second curvature}. \textit{Curvature}\footnote{The notion of \textit{curvature} of a plane curve first appears in the work of Apollonius of Perga.} measures, so to speak,
the deviation of the curve from a straight line in the neighborhood
of any of its points. While the \textit{torsion}\footnote{The name \textit{torsion} is due to L.I. Vallée, \textit{Traité de Géométrie Descriptive}.} measures, roughly
speaking, the magnitude and sense of deviation of the curve from the
\textit{osculating plane}\footnote{ The \textit{osculating plane} is
defined as the plane spanned by the instantaneous velocity and
acceleration vectors.} in the neighborhood of the corresponding
point of the curve, or, in other words, the rate of change of the
\textit{osculating plane}. Physically, a three-dimensional curve may
be obtained from a straight line by bending (\textit{curvature}) and
twisting (torsion). For high dimensions greater than three, say $n$,
a $n$-dimensional curve has $\left( {n - 1} \right)$
\textit{curvatures} which may be computed while using the
Gram-Schmidt orthogonalization process \cite{gluck} and provides the
Fr\'{e}net formulas \cite{fren} for a $n$-dimensional curve.\\

In \cite{gin1} it has been established that the location of the point where the \textit{curvature of the flow}, i.e. the \textit{curvature} of the \textit{trajectory curves }integral of any \textit{slow-fast dynamical systems} of low dimensions two and three vanishes directly provides the \textit{slow} \textit{invariant } \textit{manifold} analytical equation associated to such dynamical systems.\\

In a book recently published \cite{gin3} the \textit{Flow Curvature Method} has been generalized to high-dimensional dynamical systems and then extensively exemplified to analytically compute: \textit{fixed points stability}, \textit{invariant sets}, \textit{center manifold approximation}, \textit{normal forms}, \textit{local bifurcations}, \textit{linear invariant manifolds} of any $n$-\textit{dimensional dynamical systems} which may be used to build \textit{first integrals} of these systems.\\

One of the main applications of the \textit{Flow Curvature Method} presented in the next section establishes that \textit{curvature of the flow}, i.e. \textit{curvature} of \textit{trajectory curves} of any $n$-dimensional dynamical system directly provides its \textit{slow
manifold} analytical equation the \textit{invariance} of which is
proved according to \textit{Darboux Theorem}.\\

Then, it will be used for the first time to analytically compute the \textit{slow invariant manifold analytical equation} of the four-dimensional Unforced and Forced Heartbeat Model.

\section{Slow Invariant Manifold Analytical Equation}

The concept of \textit{invariant manifolds} plays a very important role in the stability and structure
of dynamical systems and especially for \textit{slow-fast dynamical systems} or \textit{singularly perturbed systems}. Since the beginning of the
twentieth century it has been subject to a wide range of seminal research.
The classical geometric theory developed originally by Andronov \cite{an},
Tikhonov \cite{tik} and Levinson \cite{lev} stated that \textit{singularly perturbed systems} possess \textit{invariant manifolds} on which
trajectories evolve slowly and toward which nearby orbits contract
exponentially in time (either forward and backward) in the normal
directions. These manifolds have been called asymptotically stable (or
unstable) \textit{slow manifolds}. Then, Fenichel \cite{fen1}, \cite{fen4} theory for the persistence of normally hyperbolic
invariant manifolds enabled to establish the local invariance of \textit{slow manifolds} that
possess both expanding and contracting directions and which were labeled
\textit{slow invariant manifolds}.\\
Thus, various methods have been developed in order to determine the \textit{slow invariant}
\textit{manifold} analytical equation associated to \textit{singularly perturbed systems}. The essential works of Wasow \cite{wasow},
Cole \cite{cole}, O'Malley \cite{malley1}, \cite{malley2} and Fenichel \cite{fen1}, \cite{fen4} to name but a
few, gave rise to the so-called \textit{Geometric Singular Perturbation Theory} and the problem for finding the \textit{slow invariant manifold} analytical
equation turned into a regular perturbation problem in which one generally
expected, according to O'Malley (1974 p. 78, 1991 p. 21) the asymptotic
validity of such expansion to breakdown.\\
So, the main result of this work established in the next section is that
\textit{curvature of the flow}, i.e. \textit{curvature} of
\textit{trajectory curves} of any $n$-dimensional dynamical system
directly provides its \textit{slow manifold} analytical equation the
\textit{invariance} of which is established according to
\textit{Darboux Theorem}. Since it uses neither eigenvectors nor
asymptotic expansions but simply involves time derivatives of the
velocity vector field, it constitutes a general method simplifying
and improving the \textit{slow invariant manifold} analytical
equation determination of high-dimensional dynamical systems.

\subsection{Slow manifold of high-dimensional dynamical systems}

In the framework of \textit{Differential Geometry} \textit{trajectory curves} $\vec
{X}\left( t \right)$ integral of $n$-dimensional dynamical systems (\ref{eq1})
satisfying the assumptions of the Cauchy-Lipschitz theorem may be regarded
as $n$-dimensional \textit{smooth curves}, i.e. \textit{smooth curves} in Euclidean $n-$space \textit{parametrized in terms of time}.\\

\begin{prop}
\label{prop1}
\textit{The location of the points where the curvature of the flow,
i.e. the curvature of the trajectory curves of any n-dimensional
dynamical system vanishes directly provides its $\left( {n - 1}
\right)$-dimensional slow invariant manifold analytical
equation which reads:}

\begin{eqnarray}
\label{eq2}
\phi ( \vec {X} ) & = & \dot {\vec {X}} \cdot
\left( {\ddot {\vec {X}} \wedge \dddot {\vec {X}} \wedge \ldots
\wedge \mathop {\vec {X}}\limits^{\left( n \right)} } \right) \\ \nonumber
& = & det \left( { \dot {\vec {X}},\ddot {\vec {X}},\dddot {\vec {X}},\ldots,\mathop {\vec {X}}\limits^{\left( n \right)} } \right) = 0
\end{eqnarray}

\textit{where} $\mathop {\vec {X}}\limits^{\left( n \right)}$ \textit{represents the time derivatives of} ${\vec {X}}$.

\end{prop}

\vspace{0.1in}

\begin{proof}
Cf. Ginoux \textit{et al.} \cite{gin2} ; Ginoux \cite{gin3}
\end{proof}

\subsection{Darboux invariance theorem}

According to Schlomiuk \cite{schlo1}, \cite{schlo2} and Llibre \textit{et al.} \cite{llibre} it seems that in his memoir
entitled: \textit{Sur les \'{e}quations diff\'{e}rentielles alg\'{e}briques du premier ordre et du premier degr\'{e},} Gaston Darboux (1878, p. 71) has been the
first to define the concept of \textit{invariant manifold}. Let's consider a $n$-dimensional dynamical
system (\ref{eq1}) describing ``the motion of a variable point in a space of
dimension $n$.'' Let $\vec {X} = \left[ {x_1 ,x_2 ,\ldots ,x_n } \right]^t$ be
the coordinates of this point and $\overrightarrow V = \left[ {\dot {x}_1
,\dot {x}_2 ,\ldots ,\dot {x}_n } \right]^t$ its velocity vector.\\

\begin{prop}
\label{prop2}
\textit{Consider the manifold defined by $\phi ( \vec {X} ) = 0$
where $\phi $ is a $C^1$ in an open set U is invariant with
respect to the flow of (\ref{eq1}) if there exists a $C^1$ function
denoted $K( \vec {X} )$ and called cofactor which
satisfies:}

\begin{equation}
\label{eq3}
L_{\overrightarrow V } \phi ( \vec {X} ) = K ( \vec {X} ) \phi ( \vec {X} )
\end{equation}

\textit{for all $\vec {X} \in U$ and with the Lie derivative operator defined as:}

\[
L_{\overrightarrow V } \phi = \overrightarrow V \cdot
\overrightarrow \nabla \phi = \sum\limits_{i = 1}^n {\frac{\partial
\phi }{\partial x_i }\dot {x}_i } = \frac{d\phi }{dt}.
\]

In the following \textit{invariance} of the \textit{slow manifold} will be established according to what will be
referred as \textit{Darboux Invariance Theorem}.\\
\end{prop}

\begin{proof}
Cf. Ginoux \textit{et al.} \cite{gin2} ; Ginoux \cite{gin3}
\end{proof}

\newpage

\section{Heartbeat model}

\subsection{Description of the model}

According to di Bernardo \textit{et al.} \cite{di} "The cardiac conduction system may be assumed to be a network of self-excitatory pacemakers, with the SinoAtrial (SA) node having the highest intrinsic rate. Subsidiary pacemakers with slower firing frequencies are located in the AtrioVentricular (AV) node and the His-Purkinje system. Under physiological conditions, the SA node is the dominant pace-maker and impulses travel from this node to the ventricule through the AV junction, which is traditionally regarded as a passive conduit." Then, starting from the assumptions (\cite{guev}) that between the SA and AV node a bi-directional coupling exists they describe the cardiac conduction system (Cf. Fig. 1 {\&} Fig. 2) by means of two-coupled nonlinear oscillators. For a \textit{genesis of the model} see di Bernardo \textit{et al.} \cite{di} and Signorini \textit{et al.} \cite{sign}.

\begin{figure}[htbp]
\centerline{\includegraphics[width=9.74cm,height=6.35cm]{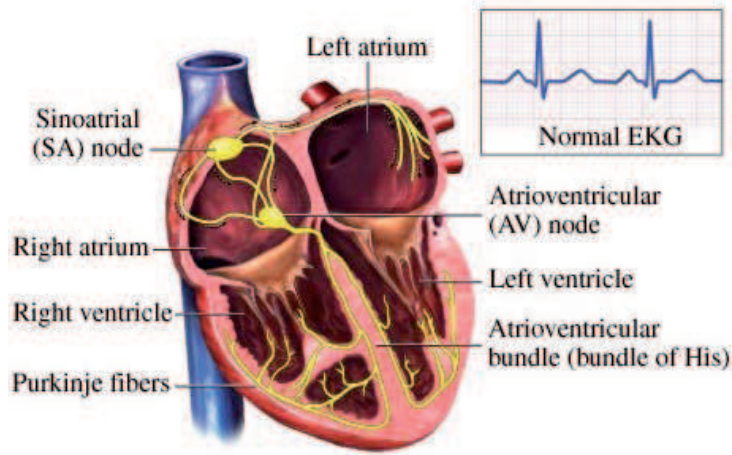}}
\caption{Cardiac Conduction System.}
\label{fig1}
\end{figure}

\begin{figure}[htbp]
\centerline{\includegraphics[width=7.35cm,height=6.27cm]{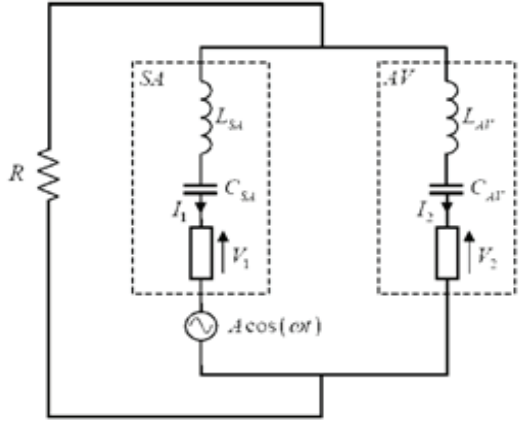}}
\caption{Equivalent Electrical Circuit.}
\label{fig2}
\end{figure}

The model proposed may be also described by an equivalent electrical circuit depicted in Fig. 2. Its structure with two nonlinear oscillators in parallel with a resistance put together the features of the Van der Pol \cite{VdP26} and West \cite{west} models. These two coupled nonlinear oscillators are built from a modification of the Van der Pol model, so that the generated waveforms resemble the action potentials of cells in the SA (resp. AV) node. The AV nonlinear oscillator represents the waveform of the intensity ($x_3$) which satisfactorily replicate the action potential ($x_4$) in the AV node and corresponds exactly to the original Van der Pol model when the resistance is null, i.e. in the uncoupled case. Thus, the nonlinear function $f(x)$ which represents the tension $(x_3)$-current $(x_4)$ characteristic of the nonlinear resistor (\textit{e.g. a triode}) may be written as:

\[
f(x) = x^3/3 - x
\]

For the generated waveform of the intensity ($x_2$) resembles the action potential ($x_1$) of the SA node cells the tension $(x_1)$-current $(x_2)$ characteristic of the nonlinear resistor of the oscillator representing the SA node which is also of Van der Pol type must be modified (di Bernardo \textit{et al.} \cite{di}). Thus, by adding the piecewise linear function $h(x)$ to the cubic nonlinearity of the Van der Pol oscillator $f(x)$ its characteristics reads.

\[
g(x) = h(x) + f(x)
\quad
\mbox{where}
\quad
h\left( x \right) = \left\{ \begin{aligned}
& +x & \mbox{ for } & x < -0.5 \hfill \\
& -x^2-0.25 &\mbox{ for } & \left| {x } \right| \leqslant 0.5 \hfill \\
& -x & \mbox{ for } & x > 0.5
\end{aligned}
\right.
\]

Moreover, let suppose that the waveform of the voltage generator in (SA) can be sinusoidal of amplitude $A$ and frequency $f$. Forcing the system means that a region of the cardiac tissue can become an active pacemaker and so interferes with the normal sinus rhythm generated by the SA node. The model is thus two-coupled nonlinear oscillator implemented in a set of four non-autonomous ordinary differential equations.

\begin{eqnarray}
\label{eq4}
\dot {\vec {X}} \left( {{\begin{array}{*{20}c}
{\dot{x_1 }} \vspace{4pt} \\
{\dot{x_2 }} \vspace{4pt} \\
{\dot{x_3 }} \vspace{4pt} \\
{\dot{x_4 }}  \\
\end{array} }} \right) & = & \left( {{\begin{array}{*{20}c}
{\frac{1}{ C_{SA} } x_2} \vspace{4pt} \\
{-\frac{1}{ L_{SA} } [x_1 + g(x_2) + R(x_2 + x_4)] + A cos(2 \pi f t) } \vspace{4pt} \\
{\frac{1}{ C_{AV} } x_2} \vspace{4pt} \\
{-\frac{1}{ L_{AV} } [x_3 + f(x_4) + R(x_2 + x_4)]} \\
\end{array} }} \right)
\end{eqnarray}

\vspace{0.1in}

The $R$ parameter models the coupling "strength" between the SinoAtrial (SA) and the AtrioVentricular (AV) node. The parameters value satisfying a normal heartbeat dynamics are:

\begin{equation}
\label{eq5}
C_{SA} = 0.25 F \mbox{; } L_{SA} = 0.05 H \mbox{; } C_{AV} = 0.675 F \mbox{; } L_{AV} = 0.027 H \mbox{; } R = 0.11 \Omega
\end{equation}

By varying this coupling resistance while keeping other parameters as above, a type of arrythmia known as $2^{o}$ \textit{AV block of the Wenckebach type} may be obtained. In order to describe arrythmia, di Bernardo \textit{et al.} \cite{di} introduced two integers $n:m$ which means the atria contract $n$ times while the ventricles $m$ times. Thus, the \textit{Flow Curvature Method} will enable, according to Prop. \ref{prop1}, to directly compute the \textit{slow manifold analytical equation} associated with heartbeat model (\ref{eq4}) in both \textit{unforced} and \textit{forced} cases for various values of the coupling parameter $R$.

\subsection{Unforced Heartbeat model}

While posing $A=0$ in Eq. (\ref{eq4}) the heartbeat model is then described by an autonomous dynamical system. As soon as $R > 0.11 \Omega$, $1:1$ periodic solutions are observed (Cf. Fig. 3).

\begin{figure}[htbp]
\centerline{\includegraphics[width=7cm,height=7cm]{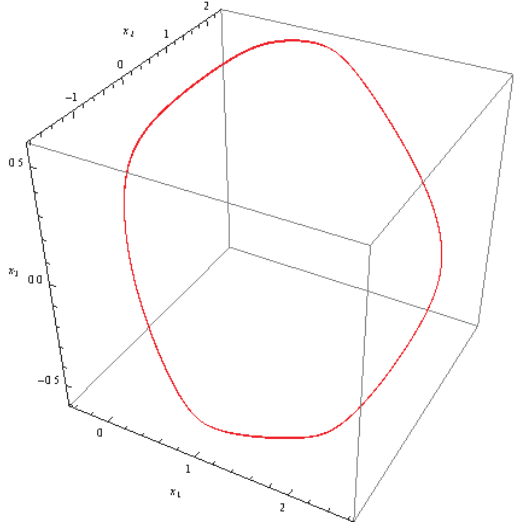}}
\caption{Unforced Heartbeat Model trajectory in the $(x_1, x_2, x_3)$ phase space.}
\label{fig3}
\end{figure}

When $R$ decreases this value a series of subharmonic bifurcations undergoes and the attractor resembles one type of arrythmia known as $2^{o}$ \textit{AV block of the Wenckebach type}. By posing $R=0.018 \Omega$ the \textit{slow invariant manifold analytical equation} corresponding to that case (Cf. Fig. 4) has been computed according to Prop. \ref{prop1} and may be written as:

\begin{equation}
\label{eq6}
\phi ( \vec {X} ) = det\left( { \dot {\vec {X}},\ddot {\vec {X}},\dddot {\vec {X}}, \ddddot {\vec {X}} } \right) = 0
\end{equation}

Because of the presence of the piecewise linear function $h(x)$ in Eq. (\ref{eq4}) the \textit{slow invariant manifold analytical equation} may be computed for each side, i.e. for $x<-0.5$ and for $x>0.5$. In both cases it leads to a polynomial depending on the variables: $x_1,x_2, x_3, x_4$ the equation of which is given in Appendix.

\begin{figure}[htbp]
\centerline{\includegraphics[width=7cm,height=7cm]{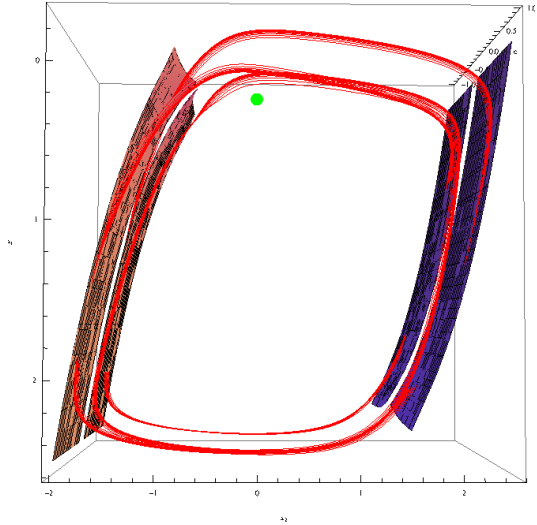}}
\caption{Slow invariant manifold of the Unforced Heartbeat Model trajectory in the $(x_1, x_2, x_3)$ phase space.}
\label{fig4}
\end{figure}

In Fig. 4 it can be observed that both trajectory curves and \textit{slow invariant manifold} are in a very close vicinity. This due to the smallness of the parameters such that $L_{SA}$ and $L_{AV}$.

\subsection{Forced Heartbeat model}

When the amplitude $A$ is different from zero the dynamical system (\ref{eq4}) will become non-autonomous and it will be the same for the \textit{slow invariant manifold analytical equation}. So, in order to avoid such difficulty a suitable variable changes may transform this \textit{non-autonomous} system into an \textit{autonomous} one while increasing the dimension of two. Indeed the forcing is modeled by the cosine which is nothing else but the solution of an harmonic oscillator. Thus, the \textit{autonomous forced heartbeat model} may be written as a set of six ordinary differential equations.

\begin{eqnarray}
\label{eq7}
\dot {\vec {X}} \left( {{\begin{array}{*{20}c}
{\dot{x_1 }} \vspace{4pt} \\
{\dot{x_2 }} \vspace{4pt} \\
{\dot{x_3 }} \vspace{4pt} \\
{\dot{x_4 }} \vspace{4pt} \\
{\dot{x_5 }} \vspace{4pt} \\
{\dot{x_6 }}  \\
\end{array} }} \right) & = & \left( {{\begin{array}{*{20}c}
{\frac{1}{ C_{SA} } x_2} \vspace{4pt} \\
{-\frac{1}{ L_{SA} } [x_1 + g(x_2) + R(x_2 + x_4)] + A x_5 } \vspace{4pt} \\
{\frac{1}{ C_{AV} } x_2} \vspace{4pt} \\
{-\frac{1}{ L_{AV} } [x_3 + f(x_4) + R(x_2 + x_4)]} \\
{\Omega x_6} \\
{-\Omega x_5} \\
\end{array} }} \right)
\end{eqnarray}

where $x_5$ is the solution of the harmonic oscillator of pulsation $\Omega = 2\pi f$.

\begin{figure}[htbp]
\centerline{\includegraphics[width=7cm,height=7cm]{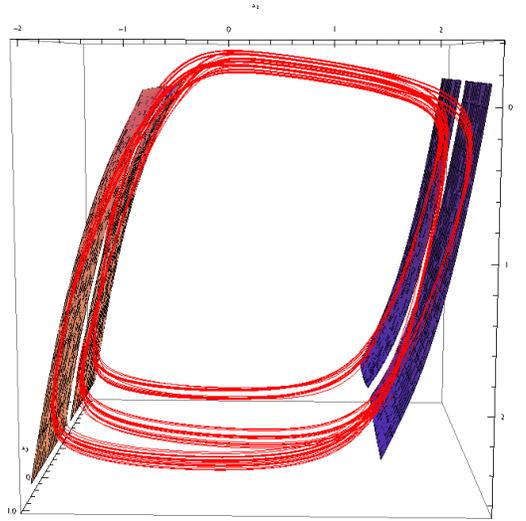}}
\caption{Slow invariant manifold of the Forced Heartbeat Model trajectory in the $(x_1, x_2, x_3)$ phase space.}
\label{fig5}
\end{figure}

Moreover, it may be stated that in the vicinity of the \textit{flow curvature manifold} both \textit{flow curvature manifold} and its Lie derivative are merged. Thus, according to \textit{Darboux Invariance Theorem} and Prop. \ref{prop2} the \textit{slow manifold} of the Unforced and Forced Heartbeat models are \textit{locally invariant}.

\section{Discussion}

In this work a new approach which consists in applying \textit{Differential Geometry} to \textit{Dynamical Systems} and called \textit{Flow Curvature Method} has been partially presented. By considering the \textit{trajectory curve}, integral of any $n$-dimensional dynamical system, as a \textit{curve} in Euclidean $n$-space, the \textit{curvature} of the \textit{trajectory curve}, i.e. \textit{curvature of the flow} has been analytically computed enabling thus to define a manifold called: \textit{flow curvature manifold}. Since such \textit{manifold} only involves the time derivatives of the velocity vector field and so, contains information about the \textit{dynamics} of the \textit{system}, it enables to find again the main features of the \textit{dynamical system} studied.  Thus, \textit{Flow Curvature Method} enables to analytically compute: \textit{fixed points stability}, \textit{invariant sets}, \textit{center manifold approximation}, \textit{normal forms}, \textit{local bifurcations}, \textit{slow invariant manifold} and \textit{integrability} of any $n$-\textit{dimensional dynamical systems} but also to ``detect'' \textit{linear invariant manifolds} of any $n$-\textit{dimensional dynamical systems} which may be used to build \textit{first integrals} of these systems.\\

Then, according to \cite{gin3} \textit{Flow Curvature Method} has been applied to a $4$-\textit{dimensional autonomous dynamical system}, i.e. the Unforced Heartbeat model and to a $6$-\textit{dimensional non-autonomous dynamical systems}, i.e. the Forced Heartbeat model in order to compute their \textit{slow invariant manifold analytical equation} which may be considered as a "state equation" linking all variables of such systems and then allowing to express one with respect to all others. Thus, according to the works of di Bernardo \textit{et al.} \cite{di}; Signorini \textit{et al.} \cite{sign} the \textit{slow invariant manifold analytical equation} could be used to study different aspects of the heartbeat dynamics such that heart prediction and control of one variable or one parameter from all others. Moreover, the \textit{Flow Curvature Method} could be also applied to biodynamical model of HIV-1. This will be the subject of another publication.

%============================ References ==============================

\newpage

\section*{Appendix}

This appendix provides the \textit{slow invariant manifold analytical equation} of the Unforced Heartbeat model for the left and right side, i.e. for $x < -0.5$ and $x > 0.5$. For sake of simplicity variables have been taken such that $x_1 = x, x_2 = y , x_3 = z, x_4 = 0$:

\[
\begin{aligned}
& \phi_{left} (x, y ,z, 0 ) = 0.0012345 x^4 - 0.0022228 x^3 y + 0.00037511 x^2 y^2 \\
& - 0.00044596 x y^3 + 0.00037448 x^3 y^3 + 0.00003387 y^4 - 0.0022672 x^2 y^4 \\
& + 0.000974419 x y^5 + 0.00033587 y^6 - 0.0066214 x^2 y^6 - 0.0033732 x y^7 \\
& + 0.00049183 y^8 - 0.0056278 x y^9 - 0.0009777 y^{10} - 0.0011415 y^{12} \\
& - 0.30904 x^3 z + 1. x^5 z + 0.31868 x^2 y z + 0.108 x^4 y z - 0.055861 x y^2 z \\
& - 3.6691 x^3 y^2 z + 0.064882 y^3 z + 0.47798 x^2 y^3 z + 4.6666 x^4 y^3 z \\
& - 0.077528 x y^4 z - 5.1675 x^3 y^4 z + 0.066041 y^5 z - 3.7898 x^2 y^5 z \\
& - 0.031419 x y^6 z + 3.1111 x^3 y^6 z - 0.040448 y^7 z - 5.3475 x^2 y^7 z \\
& - 1.4508 x y^8 z - 0.077657 y^9 z + 0.37037 x^2 y^9 z - 1.8265 x y^{10} z \\
& - 0.19842 y^{11} z - 0.16049 x y^{12} z - 0.20739 y^{13} z - 0.032921 y^{15} z \\
& + 18.784 x^2 z^2 - 102.02 x^4 z^2 + 0.35642 x y z^2 - 200.36 x^3 y z^2 \\
& + 3.7585 y^2 z^2 + 34.327 x^2 y^2 z^2 + 55.555 x^4 y^2 z^2 - 7.2815 x y^3 z^2 \\
& - 439.12 x^3 y^3 z^2 + 3.4234 y^4 z^2 - 188.62 x^2 y^4 z^2 + 5.0932 x y^5 z^2 \\
& - 37.037 x^3 y^5 z^2 - 3.9484 y^6 z^2 - 377.1 x^2 y^6 z^2 - 63.569 x y^7 z^2 \\
& - 2.1979 y^8 z^2 - 74.074 x^2 y^8 z^2 - 120.14 x y^9 z^2 - 7.6519 y^{10} z^2 \\
& - 28.806 x y^{11} z^2 - 13.151 y^{12} z^2 - 3.4293 y^{14} z^2 + 0.67914 x z^3 \\
& - 8.2025 x^3 z^3 - 0.016924 y z^3 + 32.686 x^2 y z^3 + 0.87526 x y^2 z^3 \\
& + 9.2199 y^3 z^3 - 36.513 x^2 y^3 z^3 - 26.337 x y^4 z^3 + 4.3502 y^5 z^3 \\
& - 47.011 x y^6 z^3 - 12.868 y^7 z^3 - 11.917 y^9 z^3 + 0.0061728 z^4 \\
& + 647.95 x^2 z^4 + 12.536 x y z^4 + 129.58 y^2 z^4 - 352.81 x^2 y^2 z^4 \\
& - 235.03 x y^3 z^4 + 61.791 y^4 z^4 - 588.02 x y^5 z^4 - 156.69 y^6 z^4 \\
& - 156.8 y^8 z^4 + 11.76 x z^5 - 7.8403 y^3 z^5
\end{aligned}
\]

\newpage

\[
\begin{aligned}
& \phi_{right} (x, y ,z, 0 ) = 0.0012345 x^4 - 0.004692 x^3 y + 0.0048652 x^2 y^2 \\
& - 0.00078467 x y^3 + 0.00037448 x^3 y^3 + 0.00003387 y^4 + 0.0076833 x^2 y^4 \\
& - 0.03568 x y^5 + 0.037643 y^6 - 0.0066214 x^2 y^6 + 0.03633 x y^7 \\
& - 0.046686 y^8 - 0.0056278 x y^9 + 0.013603 y^{10} - 0.0011415 y^{12} \\
& - 0.034699 x^3 z + 1. x^5 z + 0.14356 x^2 y z - 11.892 x^4 y z - 0.15531 x y^2 z \\
& + 54.378 x^3 y^2 z + 0.014991 y^3 z - 120.87 x^2 y^3 z + 4.6666 x^4 y^3 z + 131.23 x y^4 z \\
& - 37.167 x^3 y^4 z - 55.767 y^5 z + 110.55 x^2 y^5 z - 145.14 x y^6 z + 3.1111 x^3 y^6 z \\
& + 70.682 y^7 z - 17.347 x^2 y^7 z + 31.864 x y^8 z - 19.091 y^9 z + 0.37037 x^2 y^9 z \\
& - 0.93763 x y^{10} z + 0.35116 y^{11} z - 0.16049 x y^{12} z + 0.3852 y^{13} z \\
& - 0.032921 y^{15} z - 0.31296 x^2 z^2 + 9.0823 x^4 z^2 + 0.62489 x y z^2 - 95.29 x^3 y z^2 \\
& - 0.072986 y^2 z^2 + 330.69 x^2 y^2 z^2 + 55.555 x^4 y^2 z^2 - 469.16 x y^3 z^2 \\
& - 290.97 x^3 y^3 z^2 + 232.73 y^4 z^2 + 425.62 x^2 y^4 z^2 + 4.3027 x y^5 z^2 \\
& - 37.037 x^3 y^5 z^2 - 270.86 y^6 z^2 + 363.63 x^2 y^6 z^2 - 998.42 x y^7 z^2 \\
& + 838.7 y^8 z^2 - 74.074 x^2 y^8 z^2 + 340.76 x y^9 z^2 - 383.62 y^{10} z^2 \\
& - 28.806 x y^{11} z^2 + 62.294 y^{12} z^2 - 3.4293 y^{14} z^2 + 0.054748 x z^3 \\
& - 8.2025 x^3 z^3 + 0.036957 y z^3 + 36.913 x^2 y z^3 - 57.359 x y^2 z^3 \\
& + 31.972 y^3 z^3 - 36.513 x^2 y^3 z^3 + 191.33 x y^4 z^3 - 233.44 y^5 z^3 \\
& - 47.011 x y^6 z^3 + 110.02 y^7 z^3 - 11.917 y^9 z^3 + 0.0061728 z^4 \\
& - 57.678 x^2 z^4 + 115.19 x y z^4 - 11.168 y^2 z^4 - 352.81 x^2 y^2 z^4 \\
& + 2117 x y^3 z^4 - 2778.3 y^4 z^4 - 588.02 x y^5 z^4 + 1411.3 y^6 z^4 \\
& - 156.8 y^8 z^4 +  11.76 x z^5 - 7.8403 y^3 z^5
\end{aligned}
\]

The complete \textit{slow invariant manifold analytical equation} of the Unforced and Forced Heartbeat model may be downloaded at: http://ginoux.univ-tln.fr

\end{document}